\documentclass[12pt]{article}
\usepackage{amsbsy,amssymb,amsmath,amsfonts,epsfig}

\newtheorem{theorem}{Theorem}[section] 
\newtheorem{lemma}[theorem]{Lemma}     
\newtheorem{corollary}[theorem]{Corollary}
\newtheorem{proposition}[theorem]{Proposition}
\newtheorem{defn}[theorem]{Definition}

\newtheorem{rem}[theorem]{Remark}

\newtheorem{problem}[theorem]{Problem}

\def\pf{ \noindent {\bf Proof: \  }}

\newcommand{\qed}{\hfill\vrule height6pt width6pt depth0pt}

\def\endpf{\qed \medskip}

\newcommand{\Nb}{N}
\newcommand{\Mb}{M}
\newcommand{\ceq}{=}

\newcommand{\cI}{\mathcal{I}}

\newcommand{\AD}{\mathcal{A}}

\newcommand{\N}{\mathbb{N}}

\title{The SHAI property for the operators on $L^p$
 }

 \author{W.~B.~Johnson,\thanks{Supported in part by NSF
  DMS-1900612}\ \ \ %
N.~C.~Phillips,\thanks{Supported in part by NSF DMS-1501144
and by Simons Foundation Collaboration Grant for Mathematicians
\#587103}
 \ \ and G. Schechtman\thanks{Supported in part by the Israel
  Science Foundation}}

\date{7~February 2021}

\begin{document}
\maketitle

\begin{abstract}
A Banach space $X$ has the SHAI (surjective homomorphisms are injective)
property provided that for every Banach space $Y$,
every continuous surjective algebra homomorphism
from the bounded linear operators on $X$
onto the bounded linear operators on $Y$ is injective.
The main result gives a sufficient condition for $X$
to have the SHAI property.
The condition is satisfied for $L^p(0,1)$ for $1<p<\infty$,
spaces with symmetric bases that have finite cotype,
and the Schatten $p$-spaces for $1<p<\infty$.

\end{abstract}

\section{The main results}
Following Horvath \cite{Horv},
we say that a Banach space $X$ has the SHAI
(surjective homomorphisms are injective) property
provided that for every Banach space $Y$,
every surjective continuous algebra homomorphism
from the space $L (X)$ of bounded linear operators on $X$ onto $L (Y)$
is injective,
and hence by Eidelheit's \cite{Eidh} classical theorem,
$X$ is isomorphic as a Banach space to $Y$.
The continuity assumption is redundant
by an automatic continuity theorem
of B.~E.~Johnson \cite[Theorem 5.1.5]{dales}.
The spaces $\ell^p$ for $1\le p \le \infty$
are known to have the SHAI property \cite[Proposition 1.2]{Horv},
as do some other classical spaces \cite{Horv}, \cite{HorvKan},
but there are many spaces
that do not have the SHAI property \cite{Horv}.
Our research on the SHAI property was motivated by
the problem mentioned by Horvath \cite{Horv}
whether $L^p \ceq L^p(0,1)$ has the SHAI property.
A consequence of our main results, Corollary \ref{LpSHAI},
is that for $1<p<\infty$, the space $L^p$ has the SHAI property.
We do not know whether $L^1$ has the SHAI property.
The space $L^{\infty}$ does have the SHAI property
because $L^{\infty}$ is isomorphic as a Banach space
to $\ell^{\infty}$ \cite[Theorem 4.3.10]{AK}.

Before stating our theorems,
we need to review the notion
of an unconditional Schauder decomposition of a Banach space $X$.
A family $(E_{\alpha})_{\alpha \in A}$  of closed subspaces of  $X$
is called an unconditional Schauder decomposition for $X$
provided every vector $x$ in $X$ has a unique representation
$x = \sum_{\alpha \in A} x_{\alpha}$,
where the convergence is unconditional and, for each $\alpha \in A$,
the vector $x_{\alpha}$ is in $E_{\alpha}$.
Notice that by uniqueness of the representation,
$E_{\alpha} \cap E_{\beta} = \{0\}$ when $\alpha \neq \beta$,
and there are idempotents $P_{\alpha}$ on $X$
such that  $P_{\alpha} X = E_{\alpha}$
and $P_{\alpha} P_{\beta} = 0$ for $\alpha \neq \beta$.
It is known that the $P_{\alpha}$ are in $L (X)$.
Moreover, for any subset $B$ of $A$, the net
$\{ \sum_{\alpha \in F} P_{\alpha} \colon
 {\mbox{$F \subset B$ finite}}\}$
is bounded in $L (X)$
and converges strongly to an idempotent $P_B$
that has range  $\overline{\mathrm{span}}_{\alpha \in B}\  E_{\alpha}$.
The suppression constant of the decomposition is then defined to be
$\sup \{ \left\| \sum_{\alpha \in F} P_{\alpha} \right\| \colon
  {\mbox{$F \subset A$ finite}} \}$.
Note that $\| P_B \|$ is bounded by this suppression constant
for all subsets $B$ of~$A$.
In practice, this theorem is rarely used,
since typically one constructs the idempotents $P_{\alpha}$
and checks the uniform boundedness of the aforementioned nets
and verifies  the statement about the ranges
of the strong limits of the nets.
Finally, observe that  a collection $(e_{\alpha})_{\alpha \in A}$
forms an unconditional Schauder basis for $X$ if and only if
$(E_{\alpha})_{\alpha \in A}$
is an unconditional Schauder decomposition of $X$,
where $E_{\alpha} \ceq \mathbb{K} e_{\alpha}$ ($\mathbb{K}$
is the scalar field).
In the sequel, we will most often use
an unconditional Schauder decomposition $E_{\alpha}$
where each $E_{\alpha}$ is finite dimensional.
Such a decomposition is called an  unconditional FDD.
FDD stands for finite dimensional decomposition.
Schauder decomposititions and FDDs
are discussed in the monograph \cite[Section 1.g]{LT}.
Schauder bases, type/cotype theory,
and other concepts from Banach space theory that are used in this paper
are treated in the textbook \cite{AK}.

A concept that is particularly relevant for us
is that of bounded completeness.
An unconditional Schauder decomposition
$(E_{\alpha})_{\alpha \in A}$  for $X$
is said to be boundedly complete provided that
whenever $x_{\alpha} \in E_{\alpha}$ and
$\{ \left\| \sum_{\alpha \in F} x_{\alpha} \right\|_X \colon
  {\mbox{$F \subset A$ finite}} \}$
is bounded,
then the formal sum $\sum_{\alpha \in A}  x_{\alpha}$ converges in $X$,
which is the same as  saying that the net
$\{ \sum_{\alpha \in F} x_{\alpha} \colon
 {\mbox{$F \subset A$ finite}} \}$
converges.
A convenient condition that obviously guarantees bounded completeness
is that the decomposition has a disjoint lower $p$ estimate
for some $p<\infty$.
The decomposition $(E_{\alpha})_{\alpha \in A}$
is said to have a disjoint lower; respectively,
upper; $p$ estimate provided that there is $C<\infty$
so that whenever $x_1,\dots, x_n$ are finitely many vectors in $X$
such that for every $\alpha \in A$ there is at most one $i$
with $1\le i \le n$ for which $P_{\alpha} x_i \neq 0$,
we have
for $x= \sum_{i = 1}^n x_i$ the inequality
\[
\left\| \sum_{i = 1}^n  x_i \right\|
 \ge \frac{1}{C} \left( \sum_{i = 1}^n  \| x_i\|^p \right)^{1/p}; \quad
\text{respectively,} \quad
\left\| \sum_{i = 1}^n  x_i \right\|
 \le C \left( \sum_{i = 1}^n  \| x_i \|^p \right)^{1/p}.
\]
It is easy to see that the decomposition $(E_{\alpha})_{\alpha\in A}$
has a disjoint lower $p$ estimate with constant $C$
if and only if whenever $F_1,\dots,F_n$
are disjoint finite subsets of $A$ and $x$ is in $X$, then
\[
\|x\| \ge \frac{1}{C} \left( \sum_{j = 1}^n
  \left\| \sum_{\alpha\in F_j} P_{\alpha} x \right\|^p \right)^{1/p},
\]
where, as usual,
$P_{\alpha}$ is the idempotent associated with the decomposition.
Important for us is the following observation,
which is very easy to prove.
Suppose that $(E_{\alpha})_{\alpha \in A}$
is an unconditional Schauder decomposition for a subspace $X$
of a Banach space $Y$.
Assume that the idempotents $\tilde{P}_{\alpha}$ associated
with the decomposition extend to commuting idempotents ${P}_{\alpha}$
from $Y$ onto $E_{\alpha}$ and that the net
$\{ \sum_{\alpha \in F} P_{\alpha} \colon
 {\mbox{$F \subset A$ finite}} \}$
is bounded in $L (Y)$.
If $(E_{\alpha})_{\alpha \in A}$
is a boundedly complete unconditional Schauder decomposition  of $X$,
then  for each subset $B$ of $A$, the net
$\{ \sum_{\alpha \in F} P_{\alpha} \colon
 {\mbox{$F \subset B$ finite}} \}$
converges strongly in $L (Y)$  to an idempotent  $P_B$
whose range is the closed linear span of the spaces $E_{\alpha}$
for $\alpha \in B$
(which, by abuse of notation, we abbreviate to
$\overline{\mathrm{span}} \, \{ E_{\alpha} \colon \alpha \in B \}$)
and $P_B$ extends the basis projection from $X$ onto
$\overline{\mathrm{span}} \, \{E_{\alpha} \colon \alpha \in B\}$.
In particular, $X$ is complemented in $Y$.
Conversely, if $X$ is known to be complemented in $Y$,
then such extensions $P_B$ of the basis projections $\tilde{P}_B$
from $X$ onto
$\overline{\mathrm{span}} \, \{E_{\alpha} \colon \alpha \in B\}$
obviously exist even when the decomposition is not boundedly complete.
In general, to guarantee that $X$ is complemented in $Y$,
something is needed other than having commuting extensions $P_{\alpha}$
with
$\{ \sum_{\alpha \in F} P_{\alpha} \colon
 {\mbox{$F \subset A$ finite}} \}$
uniformly bounded:
consider $X=c_0$,
$Y=\ell^{\infty}$, and the unit vector basis of $c_0$.

From the definitions of type and cotype,
it is clear that if $X$ has type $p$ and cotype $q$,
then every unconditional Schauder decomposition for $X$
has a disjoint upper $p$ estimate and a disjoint lower $q$ estimate,
where the constants depend only on the suppression constant
of the decomposition and the  type $p$ and cotype $q$ constants of $X$.
In particular, if $1<p\le 2$,
then every unconditional Schauder decomposition
for a subspace of a quotient of $L^p$
has a disjoint upper $p$ estimate and a disjoint lower $2$ estimate,
while if $2\le p <\infty$,
then every uncondtional Schauder decomposition
for a subspace of a quotient of $L^p$
has a disjoint upper $2$ estimate
and a disjoint lower $p$ estimate \cite[Theorem 6.2.14]{AK}.

The observation in the following lemma
will be used for transferring information from $Y$ to $X$
when there is a surjective homomorphism from $L (Y)$ onto $L (X)$.
\begin{lemma}\label{lemma:pestimate}
Suppose that $(E_{\alpha})_{\alpha \in A}$
is an unconditional decomposition for $X$
that has a disjoint lower $p$ estimate with $1\le p <\infty$,
and let $Y\supseteq X$.
Then there is a constant $C<\infty$ such that
if $A_1,\dots, A_n$ are disjoint subsets of $A$
and $P_{A_j}$ is the basis projection onto
$E_{A_j} \ceq
 \overline{\mathrm{span}} \, \{E_{\alpha} \colon \alpha \in A_j\}$
and $T_1,\dots ,T_n$ are operators in $L (Y)$, then
%
\[
\left\| \sum_{i = 1}^n T_i P_i \right\|
 \le C \left( \sum_{i = 1}^n \|T_i \|^q \right)^{1/q},
\  \text{where} \ 1/p + 1/q =1.
\]
\end{lemma}
\pf Suppose $x\in X$.
Then
\begin{align*}
\left\| \sum_{i = 1}^n T_i P_i  x \right\|
& \le \sum_{i = 1}^n \| T_i \| \| P_i x  \|
\le
\left( \sum_{i = 1}^n \| T_i \|^q \right)^{1/q}
  \left( \sum_{i = 1}^n \|P_i x \|^p \right)^{1/p}
\\
& \le C \left( \sum_{i = 1}^n \| T_i \|^q \right)^{1/q}  \|x \|,
\end{align*}
where the constant $C$ is  the disjoint lower $p$ constant
of $(E_{\alpha})_{\alpha \in A}$.
\endpf

\noindent
A family of sets is said to be  \emph{almost disjoint}
provided the intersection of any two of them is finite.
\begin{defn}\label{def:propertysharp}
Suppose that $(E_n)_{n=1}^{\infty}$
is an unconditional FDD for a Banach space $X$.
We say that $(E_n)$ has property (\#) provided there is
an almost disjoint continuum $\{\Nb_{\alpha} \colon \alpha < c\}$
of infinite sets of natural numbers such that for each $\alpha < c$,
$X$ is isomorphic to the closed linear span of the subspaces
$E_n$ for $n \in \Nb_{\alpha}$.
\end{defn}
Subsymmetric bases are obvious examples of FDDs that have property (\#).
(A basis is subsymmetric if it is unconditional
and every subsequence of the basis is equivalent to the basis.
Symmetric bases are subsymmetric.)
A second almost obvious example
is the direct sum of two Banach spaces with subsymmetric bases.
Such a space has an FDD with property (\#)
such that each space in the decomposition is two dimensional.
In Corollary \ref{LpSHAI} we point our that the Haar basis for $L^p$
has property (\#) when $1<p<\infty$.
\begin{proposition}\label{prop0.2}
Let $(E_n)_{n=1}^{\infty}$ be an FDD for a Banach space $X$.
Assume that $(E_n)$ has property (\#),
witnessed by an almost disjoint family
$\{\Nb_{\alpha} \colon \alpha < c\}$
of infinite subsets of the natural numbers.
For $F\subset \N$,
let $P_F$ be the basis projection from $X$
onto the closed linear span $E_F$ of the subspaces $E_n$
for $n \in F$.
Suppose that $\Phi$ is a non zero,
non injective continuous homomorphism
from $L (X)$ onto a Banach algebra $\AD$.
Then for each $\alpha <c$,
$\Phi (P_{\Nb_{\alpha}})$ is a non zero idempotent in $\AD$.
Moreover, there is a constant $C<\infty$
such that if $F$ is any finite subset of $[\alpha <c]$, then
$\left\| \sum_{\alpha \in F} \Phi(P_{\Nb_{\alpha}}) \right\|_{\AD}
 \le  C$.
If $\AD$ is a subalgebra of $L (Y)$ for some Banach space $Y$,
then $(\Phi (P_{\Nb_{\alpha}}))_{\alpha < c}$
is a family of commuting extensions to $Y$
of the projections associated with
an unconditional Schauder decomposition for a subspace $Y_0$ of $Y$.
\end{proposition}
\pf
Since, for each $\alpha$, the range of $P_{\Nb_{\alpha}}$
is isomorphic to $X$, and $\Phi$ is not zero,
$\Phi (P_{\Nb_{\alpha}})$ is a non zero idempotent in $\AD$.
Suppose that $F$ is a finite subset of $\{ \alpha \colon \alpha < c \}$.
Take a finite set $S$ of natural numbers so that
$\Nb_{\alpha} \cap \Nb_{\beta} \subset S$
for all distinct $\alpha, \beta$ in $F$.
For $\alpha \in F$,
let $Q_{\alpha} \ceq P_{\Nb_{\alpha} \setminus S}$
be the basis projection from $X$ onto
$\overline{\mathrm{span}}
  \{E_n \colon n\in \Nb_{\alpha}\setminus S \}$.
The kernel of $\Phi$ is a non trivial ideal in $L (X)$
and hence contains the finite rank operators.
Since $P_{\Nb_{\alpha}} - Q_{\alpha}$ is a finite rank operator,
$\Phi(P_{\Nb_{\alpha}}) = \Phi(Q_{\alpha}) $ for each $\alpha \in F$.
But the projections $Q_{\alpha}$, for $\alpha \in F$,
are projections onto the closed spans of disjoint subsets
of the FDD $(E_n)_{n = 1}^{\infty}$,
so
\[
\left\| \sum_{\alpha \in F} \Phi(Q_{\alpha}) \right\|_{\AD}
 \le \left\| \sum_{\alpha \in F} Q_{\alpha} \right\| \|\Phi\|
 \le C \|\Phi\|,
\]
where  $C$ is the suppression constant of $(E_n)$.
The last statement is now obvious.
\endpf

With the preliminaries out of the way,
we state the main theorem in this article.
\begin{theorem}\label{SHAI}
Let $(E_n)_{n = 1}^{\infty}$ be an unconditional FDD
for a Banach space $X$.
Assume that $(E_n)_{n = 1}^{\infty}$ has property (\#)
(Definition~\ref{def:propertysharp})
and $(E_n)_{n = 1}^{\infty}$ has a disjoint lower $p$
estimate for some $p<\infty$.
Then $X$ has the SHAI property.
\end{theorem}
\pf
Suppose, for contradiction,
that $\Phi$ is a  non injective continuous homomorphism
from $L (X)$ onto $L (Y)$ for some non zero Banach space $Y$.
We continue with the set up in Proposition \ref{prop0.2},
where property (\#) for $(E_n)$
is witnessed by an almost disjoint family
$\{\Nb_{\alpha} \colon \alpha < c\}$
of infinite subsets of the natural numbers,
and for $F\subset \N$,
the basis projection from $X$ onto the closed linear span
$E_F$ of $\{E_n \colon n \in F\}$ is denoted by $P_F$.

We claim that to get a contradiction it is enough
to prove that the subspace $Y_0$ is complemented in $Y$.
Indeed, if $Y_0$ is complemented in $Y$,
then $L (Y_0)$ is isomorphic as a Banach algebra
to a subalgebra of $L (Y)$.
However, defining $Y_{\alpha} \ceq \Phi(P_{\Nb_{\alpha}})Y$
for $\alpha < c$,
we know that $(Y_{\alpha})_{\alpha < c }$
is an unconditional Schauder decomposition for $Y_0$.
But then for every subset $S$ of $\{ \alpha \colon \alpha < c \}$
there is an idempotent $Q_S$ from $Y_0$
onto $\overline{\mathrm{span}}\{Y_{\alpha} \colon \alpha \in S\}$
with $Q_S $ zero on all $Y_{\beta}$ for which $\beta \not\in S$.
Thus if $S_1$ and $S_2$ are different subsets
of $\{ \alpha \colon \alpha < c \}$,
then $\|Q_{S_1} - Q_{S_2} \| \ge 1$,
and hence the density character of $L (Y_0)$,
whence also of $L (Y)$, is at least $2^c$.
However, since $X$ is separable,
the density character of  $L (X) $ is at most $c$
(actually, equal to $c$ since $X$  has an uncondtional FDD),
so $L (Y)$ cannot be a continuous image of $L (X)$.
This completes the proof of the claim.

To show that $Y_0$ must be complemented in $Y$,
we use the fact proved in Proposition \ref{prop0.2}
that there is a constant $C$ such that for every finite subset
$F$ of $\{ \alpha \colon \alpha < c \}$ we have
$\left\| \sum_{\alpha \in F} \Phi(P_{\Nb_{\alpha}}) \right\|_{L (Y)}
  \le  C$.
It was remarked in the introduction to this section
that this condition guarantees that $Y_0$ is complemented in $Y$
when $(Y_{\alpha})_{\alpha<c}$ is a boundedly complete decomposition.
To see that $(Y_{\alpha})_{\alpha<c}$  is boundedly complete,
we use Lemma \ref{lemma:pestimate}.
We guarantee bounded completeness by proving that
$(Y_{\alpha})_{\alpha < c }$ has a disjoint lower $p$ estimate.
That is, we just need to find a constant $C$ so that if
$F_1, \dots, F_m$ are disjoint finite subsets
of $\{ \alpha \colon \alpha < c \}$
and $y$ is in $Y$ (or even just in $Y_0$),
then
\begin{equation}\label{lowerpestimate2}
\|y\| \ge \frac{1}{C} \left( \sum_{j = 1}^m
 \left\| \sum_{\alpha \in F_j} \Phi(P_{\Nb_{\alpha}})  y \right\|^p
     \right)^{1/p}.
\end{equation}
Just as in the proof Proposition \ref{prop0.2}, we can write
$\sum_{\alpha \in F_j}\Phi(P_{\Nb_{\alpha}})  = \Phi(Q_j)$ with $Q_j$,
for $1 \le j \le m$,
being the basis projections onto the closed spans
of disjoint sets of FDD basis spaces $(E_n)$.
So (\ref{lowerpestimate2}) can be rewritten as
\begin{equation}\label{lowerpestimate3}
\|y\|
 \ge \frac{1}{C}
   \left( \sum_{j = 1}^m \| \Phi(Q_j)  y \|^p \right)^{1/p}.
\end{equation}
From Lemma \ref{lemma:pestimate} and the surjectivity of $\Phi$,
for any $T_1, \dots, T_m$ in $L (Y)$ we have
\begin{equation}\label{lowerpestimate4}
\left\| \sum_{j = 1}^m T_j \Phi(Q_j) \right\|
 \le C \left( \sum_{j = 1}^m \|T_j\|^q \right)^{1/q},
\end{equation}
where $C$ depends only on $p$
and on $\| \Phi\|\cdot \|( \Phi^*)^{-1}\|$,
and $1/p+1/q=1$.
Take any $y\in Y$ and take $\beta_j \ge 0$ with
\[
\sum_{j = 1}^m \beta_j^q = 1
\qquad {\mbox{and}} \qquad
\sum_{j = 1}^m \beta_j \| \Phi (Q_j) y\|
 = \left( \sum_{j = 1}^m \| \Phi (Q_j) y \|^p) \right)^{1/p}.
\]
Let $y_0$ be any unit vector in $Y$ and let $T_j$ be $\Phi(Q_j)$
followed by a norm (at most) one projection onto
the (at most) one dimensional space $\mathbb{K} \Phi(Q_j) y$
followed by $\Phi(Q_j) y \mapsto \beta_j \|\Phi(Q_j)y\| y_0$.
Then by (\ref{lowerpestimate4}),
\begin{align*}
\left( \sum_{j = 1}^m  \| \Phi(Q_j)  y\|^p \right)^{1/p}
& =  \sum_{j = 1}^m  \beta_j  \| \Phi(Q_j)  y\|
  = \left\| \sum_{j = 1}^m T_j \Phi(Q_j) y \right\|
\\
& \le C \left( \sum_{j = 1}^m \|T_j \|^q \right)^{1/q} \|y\| = C \|y\|,
\end{align*}
which is (\ref{lowerpestimate3}).
\endpf

Our first corollary of Theorem \ref{SHAI} is immediate.
Its hypothesis are satisfied by many spaces that are used in analysis,
including most Orlicz and Lorentz sequence spaces.
\begin{corollary}\label{symmetricSHAI}
If $X$ has a subsymmetric basis and has finite cotype,
then $X$ has SHAI.
\end{corollary}
The next corollary solves the problem
that motivated our research into the SHAI property.
\begin{corollary}\label{LpSHAI}
For $1<p<\infty$, the space $L^p$ has the SHAI property.
\end{corollary}
\pf
In view of Theorem \ref{SHAI},
it is enough to prove that the Haar basis for $L^p$ has property (\#).
Let $\{\Nb_{\alpha} \colon \alpha < c\}$
be a continuum of almost disjoint infinite subsets
of the natural numbers $\N$.
Define for $\alpha < c$
\[
X_{\alpha}
 \ceq \overline{\mathrm{span}}
  \bigl\{ h_{n,i} \colon
    {\mbox{$n \in  \Nb_{\alpha}$ and $1 \le  i \le 2^n$}} \bigr\},
\]
where
$\{h_{n,i} \colon {\mbox{$n = 0, 1, \dots$ and $1 \le  i \le 2^n$}} \}$
is the usual (unconditional) Haar basis for $L^p(0,1)$,
indexed in its usual way,
so that $ \{|h_{n,i}| \colon 1\le i \le 2^n\}$
is the set of indicator functions
of the dyadic subintervals of $(0,1)$ that have length $2^{-n}$.
By the Gamlen--Gaudet theorem \cite{GG},
$X_{\alpha}$ is isomorphic to $L^p$
with the isomorphism constant depending only on $p$.
\endpf
\begin{rem}\label{rem:minimal}
\emph{Although our proof that $L^p$ has the SHAI property
is simple enough, it is strange.
The ``natural" way of proving that a space $X$ has the SHAI property
is to verify that for any non trivial closed ideal $\cI$ in $L (X)$,
the quotient algebra $L (X)/\cI$ contains no minimal idempotents.
(An idempotent $P$ is called minimal provided $P\neq 0$
and the only idempotents $Q$ for which $PQ=QP=Q$ are $P$ and $0$.
Rank one idempotents in $L (X)$ are minimal.)
This suggests the following problem,
which is related to the known problem
whether every infinite dimensional complemented subspace of $L^p$
is isomorphic to its square.}
\end{rem}
\begin{problem}\label{prob:minimal}
{Is there a non trivial closed ideal $\cI$ in $L (L^p)$
for which $L (L^p)/\cI$ has a minimal idempotent?}
\end{problem}

If there is a positive answer to Problem \ref{prob:minimal},
the witnessing ideal $\cI$ cannot be contained
in the ideal of strictly singular operators.
This is because every infinite dimensional complemented subspace
of $L^p$ contains a complemented subspace
that is isomorphic either to $\ell^p$ or to $\ell^2$ \cite{kp},
and the fact that idempotents in $L (X)/ \cI$
lift to idempotents in $L (X)$
when $\cI$ is an ideal that is contained in $L (X)$ \cite{barnes}.

\begin{problem}\label{prob:L1SHAI}
{Does $L^1$ have the SHAI property?}
\end{problem}

\section{Examples and permanence properties}

Here we present some more examples of spaces with property (\#)
and with the SHAI property.
We do not know whether  every complemented subspace of $L^p$
has the SHAI property,
but we show  that at least some of the known examples of such spaces do.
Along the way we state and prove some permanence properties of (\#).

The classical complemented subspaces of $L^p$
have the SHAI property when $1< p<\infty$.
This was known for $\ell^2$ and $\ell^p$ and proved above for $L^p$.
The case of $\ell^p\oplus\ell^2$
follows easily from Theorem  \ref{SHAI}.
That the remaining classical complemented subspace of $L^p$,
$\ell^p(\ell^2)$, the $\ell^p$ sum of $\ell^2$,
has (\#) and the SHAI property
follows from Proposition \ref{prop:ell^psum} below.
Before stating Proposition \ref{prop:ell^psum}
we introduce a quantitative version of property (\#).

\begin{defn}\label{def:propertyKsharp}
Suppose that  $(E_n)_{n=1}^{\infty}$
is an unconditional FDD for a Banach space $X$
and $K$ is a positive constant.
We say that $(E_n)$ has property (\#) with constant $K$
provided there is  an almost disjoint continuum
$\{\Nb_{\alpha} \colon \alpha < c\}$
of infinite sets of natural numbers such that for each $\alpha < c$,
$X$ is $K$-isomorphic to the closed linear span of
$\{E_n \colon n \in \Nb_{\alpha} \}$.
\end{defn}

Note that if $(E_n)_{n=1}^{\infty}$ has property (\#)
then it has property (\#) for some positive constant $K$.
Nevertheless, we need this quantitative notion
for the full generality of Proposition  \ref{prop:ell^psum}.

Recall that if $(e_i)$ is an unconditional basis
for some Banach space $Y$ and
$X_i$, for $i=1,2,\dots$, is a Banach space,
$\left( \bigoplus_{i = 1}^{\infty} X_i \right)_Y$
is the space of sequences $\bar x=(x_1,x_2,\dots)$ whose norm,
$\|\bar x \|
 = \big\| \sum_{i = 1}^{\infty} \|x_i \| \cdot e_i \big\|_Y$,
is finite.
We denote the subspace of
$\left( \bigoplus_{i = 1}^{\infty} X_i \right)_Y$
of all sequences of the form $(0,\dots,0,x_i,0,\dots)$
by $X_i\otimes e_i$.

\begin{proposition}\label{prop:ell^psum}
For $i = 1, 2, \dots$ let $(E_n^i)_{n=1}^{\infty}$
be an unconditional FDD for a Banach space $X_i$,
all satisfying property (\#) with a common $K$.
Then for each subsymmetric basis $(e_i)$ of some Banach space $Y$,
the unconditional FDD $(E_n^i\otimes e_i)_{i,n=1}^{\infty}$
of $\left( \bigoplus_{i = 1}^{\infty} X_i \right)_Y$ satisfies (\#).
If, in addition, the decompositions $(E_n^i)_{n=1}^{\infty}$
have  disjoint lower $p$ estimates with uniform constant
and $(e_i)$ also has such an estimate,
then $\left( \bigoplus_{i = 1}^{\infty} X_i \right)_Y$
has the SHAI property.
\end{proposition}

\pf
For each $i$,
let $\{\Nb_{\alpha}^i \colon \alpha < c\}$
be an almost disjoint continuum of infinite sets of natural numbers
such that for every $\alpha < c$,
$X_{\alpha}$ is $K$-isomorphic to the closed linear span of
the subspaces $E_n^i$ for $n \in \Nb_{\alpha}$.
Also, let $\{\Nb_{\alpha} \colon \alpha < c\}$
be an almost disjoint continuum of infinite sets of natural numbers.
Then
\[
\bigl\{ (i, n) \colon
  {\mbox{$i \in \Nb_{\alpha}$ and $n \in \Nb_{\alpha}^i$}} \bigr\}
\]
is a continuum of almost disjoint subsets of $\N\times\N$.
It is easy to see that this continuum satisfies
what is required of the unconditional FDD
$(E_n^i\otimes e_i)_{i,n=1}^{\infty}$ to satisfy (\#).
If the decompositions $(E_n^i)_{n=1}^{\infty}$
have  disjoint lower $p$ estimates with uniform constant
and $(e_i)$ also has such an estimate,
then the FDD $(E_n^i\otimes e_i)_{i,n=1}^{\infty}$
clearly has a  disjoint lower $p$ estimate as well,
so the SHAI property follows from Theorem \ref{SHAI}.
\endpf
\begin{rem}
\emph{Note that the proof above works with only notational differences
if we deal with only finitely many $X_i$
(and here we do not need to assume the uniformity of the (\#) property).
In particular,
if each of $X$ and $Y$ has an unconditional FDD with (\#),
then so does $X\oplus Y$.}
\end{rem}
As we said above, this takes care of the space $\ell^p(\ell^2)$.
The first non classical complemented subspace of $L^p$
is the space $X_p$ of Rosenthal \cite{ro}.
We recall its definition.
Let $p>2$ and let $\bar w=(w_i)_{i = 1}^{\infty}$
be a bounded sequence of positive real numbers.
Let $(e_i)_{i = 1}^{\infty}$ and $(f_i)_{i = 1}^{\infty}$
be the unit vector bases of $\ell^p$ and $\ell^2$.
Let $X_{p,\bar w}$ be the closed span of
$(e_i\oplus w_i f_i)_{i = 1}^{\infty}$ in $\ell^p\oplus\ell^2$.
If the $w_i$ are bounded away from zero,
then $X_{p,\bar w}$ is isomorphic to $\ell^2$.
If $\sum_{i = 1}^{\infty} w_i^{\frac{2p}{p-2}}<\infty$,
then  $X_{p,\bar w}$ is isomorphic to $\ell^p$.
If one can split the sequence $\bar w$ into two subsequences,
one bounded away from zero
and the other such that the sum of the $\frac{2p}{p-2}$ powers
of its elements converges,
then $X_{p,\bar w}$ is isomorphic to $\ell^p\oplus\ell^2$.
Rosenthal proved that in all other situations one gets a new space,
isomorphically unique (i.e., any,
two spaces corresponding to two
choices of $\bar w$ with this condition are isomorphic).
Moreover,
$X_{p,w}$ is isomorphic to a complemented subspace of $L^p$.
The constants involved (isomorphisms and complementations)
are bounded by a constant depending only on $p$.
This common (class of) space(s) is denoted by $X_p$.
For $1<p<2$, $X_p$ is defined to be $X_{p/(p-1)}^*$.

\begin{proposition}\label{claim:X_p}
Let $p \in (1, \infty) \setminus \{ 2 \}$.
Then $X_p$ has (\#) and has the SHAI property.
\end{proposition}

\pf
Let $p>2$.
Write $\N$ as a disjoint union of finite subsets $\sigma_j$
for $j=1,2,\dots$,
with $|\sigma_j|\to\infty$.
For $i\in \sigma_j$ put $w_i=|\sigma_j|^{\frac{2-p}{2p}}$,
so  $w_i\to 0$ and for each $j$,
$\sum_{i\in \sigma_j}w_i^{\frac{2p}{p-2}}=1$.
Set $E_j \ceq {\mathrm{span}}\,(e_i\oplus w_if_i)_{i\in \sigma_j}$.
It follows that for any infinite subsequence
of the unconditional FDD $(E_j)$,
the closed span of this subsequence is isomorphic to $X_p$.
The FDD is unconditional and,
as it lives in $L^p$, has a lower $p$ estimate.
So the result in this case follows from Theorem \ref{SHAI}.
The case $1<p<2$ follows by looking at the dual FDD.
\endpf

Building on $X_p$ and the classical complemented subspaces of $L^p$,
Rosenthal \cite{ro} lists a few more isomorphically distinct spaces
that are isomorphic to complemented subspaces of $L^p$
when $p \in (1, \infty) \setminus \{ 2 \}$.
Using the discussion above
one can easily show that they all have (\#) and the SHAI property.
Here we just  comment on one of them
for which the full power of Proposition \ref{prop:ell^psum} is needed.
This is the space denoted in \cite{ro} by $B_p$.
It is the $\ell^p$ sum of spaces $X_i$
each having a $1$-symmetric basis,
and thus having (\#) with uniform constant.
Each $X_i$ is isomorphic to  $\ell^2$,
but the isomorphism constant tends to infinity as $i \to \infty$.
By Proposition \ref{prop:ell^psum},
$B_p$  has (\#) and  the SHAI property.

The first infinite collection
of mutually non isomorphic complemented subspaces
of $L^p$ for $p \in (1, \infty) \setminus \{ 2 \}$
was constructed in \cite{sc}.
We recall the simple construction.
Given two subspaces $X$ and $Y$ of $L^p(\Omega)$
with $1\le p\le \infty$,
$X\otimes_p Y$ denotes the subspace of $L^p(\Omega^2)$
that  is the closed span of all functions
of the form $h(s,t)=f(s)g(t)$ with $f\in X$ and $g\in Y$.
It is easy to see (and was done in \cite{sc})
that the isomorphism class of $X\otimes_p Y$
depends only on the isomorphism classes of $X$ and $Y$ and that,
if $X$ and $Y$ are complemented in $L^p(\Omega)$,
then $X\otimes_p Y$ is complemented in $L^p(\Omega^2)$.
More generally,
if $X_1,X_2,Y_1,Y_2$ are subspaces of $L^p(\Omega)$
and $T_i\in L (X_i,Y_i)$,
then $T_1\otimes_p T_2\in L (X_1\otimes_p X_2,Y_1\otimes_p Y_2)$.
Note also that if $(E_n^i)_{n=1}^{\infty}$
is an unconditional FDD for $X_i$ for $i=1,2$,
then $(E_n^1\otimes_pE_m^2)_{n,m=1}^{\infty}$
is an unconditional FDD for $X_1\otimes_p X_2$.
This follows from iterating Khinchine's inequality.

With a little abuse of notation we denote by $X_p$
some isomorph of $X_p$ that is complemented in $L^p[0,1]$.
Set $Y_1=X_p$, and for $n=2,3,\dots$, let $Y_n=Y_{n-1}\otimes_p X_p$.
From the above it is clear that the spaces $Y_n$ are complemented
(alas, with norm of projection depending on $n$)
in some $L^p$ space isometric to $L^p[0,1]$.
The main point in \cite{sc}
was to prove that the spaces $Y_n$ are isomorphically different.
That all the spaces $Y_n$ have (\#)
follows now from the following general proposition,
because it is clear that $\otimes_p$ satisfies Conditions
(\ref{Item_251}) and (\ref{Item_252}) in Proposition \ref{prop:tensor}
for the class of all $m$ tuples of subspaces of $L^p(\mu)$ spaces.
\begin{proposition}\label{prop:tensor}
Assume that $X_1,\dots,X_m$ are Banach spaces,
each of which has an unconditional FDD satisfying (\#).
Let $Y_1\otimes\dots\otimes Y_m$ denote an $m$ fold tensor product
endowed with norm defined on some class of $m$ tuples
of Banach spaces with the following two properties:
\begin{enumerate}
\item\label{Item_251}
If $T_i\in L (Y_i,Z_i)$ for $i=1,\dots,m$, then
\[
T_1\otimes\dots\otimes T_m \colon
  Y_1\otimes\dots\otimes Y_m\to Z_1\otimes\dots\otimes Z_m
\]
is bounded.
\item\label{Item_252}
If $Y_i$ has an unconditional FDD $(F_n^i)_{n=1}^{\infty}$
for each~$i$,
then
$(F_{n_1}^1\otimes\dots\otimes F_{n_m}^m)_{n_1,\dots,n_m=1}^{\infty}$
is an unconditional FDD
for the completion of $Y_1\otimes\dots \otimes_m Y_m$.
\end{enumerate}
Then, if we assume in addition that $(X_1, \dots, X_m)$
is in this class,
the completion of $X_1\otimes\dots\otimes X_m$
has an unconditional FDD with (\#).
\end{proposition}

\pf
For each $i=1,\dots,m$,
let $(E_n^i)_{n=1}^{\infty}$ be an unconditional FDD
for a Banach space $X_i$
such that there is an almost disjoint continuum
$\{\Nb_{\alpha}^i \colon \alpha < c\}$ of infinite sets of $\N$
such that for each $\alpha < c$,
$X_i$ is isomorphic to the closed linear span
of the spaces $E_n^i$ for $n \in \Nb_{\alpha}^i$.

Consider the continuum
\[
\{\Nb_{\alpha}^1\times\dots\times \Nb_{\alpha}^m \colon \alpha< c\}
\]
of subsets of $\N^m$.
This is an almost disjoint family whose cardinality is the continuum.
Property~(\ref{Item_252})
of the tensor norms we consider guarantees that
$(E_{n_1}^1\otimes\dots\otimes E_{n_m}^m)_{n_1,\dots,n_m=1}^{\infty}$
is an unconditional FDD
for the completion of $X_1\otimes\dots\otimes X_m$.
Property~(\ref{Item_251}) implies that for each $\alpha<c$,
the closed linear span of
\[
(E_{n_1}^1 \otimes \dots \otimes E_{n_m}^m)_{(n_1, \dots, n_m)
      \in \Nb_{\alpha}^1 \times \dots \times \Nb_{\alpha}^m}
\]
is isomorphic to the completion of $X_1\otimes\dots\otimes X_m$.
\endpf
\begin{rem}
\emph{Note that in general Property~(\ref{Item_251})
does not imply Property~(\ref{Item_252}).
The Schatten classes $C_p$
for $p\neq 2$ are examples of tensor norms
that satisfy (\ref{Item_251}) but not (\ref{Item_252}).}
\end{rem}
We note that it is clear from Proposition \ref{prop:tensor}
that if $X_1,\dots,X_m$ are subspaces of $L^p$ for $1\le p<2$
that have  (sub)symmetric bases,
then $X_1\otimes_p\dots\otimes_p X_m$ has (\#) and  the SHAI property.
The class of subspaces of $L^p$ for $1\le p<2$
that have a symmetric basis
(i.e., the norm of a vector is invariant, up to a constant,
under all permutations and changes of signs of its coefficients)
is a rich family.
(For $p>2$,
up to isomorphism it includes only $\ell^p$ and $\ell^2$.)
Thus the class of tensor products above includes,
for example,
$\ell_{p_1} (\ell_{p_2} (\dots (\ell_{p_m}) \dots) )$
whenever $p\le p_1< p_2<\dots< p_m\le 2$.
\begin{problem}\label{Prob:BRSspacesSHAI}
Suppose $p \in (1, \infty) \setminus \{ 2 \}$
and let $X$ be a complemented subspace of $L^p$.
Does $X$ have the SHAI property?
What if, in addition, $X$ has an unconditional basis?
What if, in addition, $X$ is one of the $\aleph_1$ spaces
constructed in \cite{BRS}?
\end{problem}
We complete this section with a discussion of another class
of classical Banach spaces that have property $(\#)$
and thus also the SHAI property; namely,
the Schatten ideals $C_p$ of compact operators $T$ on $\ell^2$
for which the eigenvalues of $(T^*T)^{1/2}$ are $p$-summable.
We treat the case $1<p  <\infty$
but remark afterwards how one can prove that $C_1$
(trace class operators on $\ell^2$)  has the SHAI property.
Neither   $C_1$ nor its predual $C_{\infty}$
(compact operators on $\ell^2$)
has an unconditional FDD \cite{kw-p}
and hence these spaces do not have property $(\#)$.
In the sequel we also assume $p\neq 2$ because $C_2$,
being isometrically isomorphic to $\ell^2$, has already been discussed.

First, consider the subspace $T_p$ of $C_p$
consisting of the lower triangular matrices in $C_p$.
Here we include $p=1$ and $p=\infty$ but exclude $p=2$.
Neither $T_p$ nor $C_p$ has an unconditional basis \cite{kw-p},
but $T_p$ has an obvious unconditional FDD $(E_n)$; namely,
$E_n \ceq \mathrm{span}_{1\le j \le n} e_n\otimes e_j$; that is,
a matrix is in $E_n$ if and only if the only non zero terms
are in the first $n$ entries of the $n$-th row.
Since multiplying all entries in a row
by the same scalar of magnitude one is an isometry on $C_p$,
$(E_n)$ is even $1$-unconditional.
If $\Mb$ is an infinite subset of $\N$,
let $T_p (\Mb )$ be the closed span in $T_p$ of $(E_n)_{n \in \Mb}$.
Since $(E_n)$ is $1$-uncondtional,
$T_p (\Mb )$ is norm one complemented in $T_p$ and, similarly,
$T_p$ is isometric to a norm one complemented subspace of $T_p (\Mb )$.
The space $T_p$ is isomorphic to $\ell^p(T_p)$ \cite[p. 85]{al},
so the decomposition method \cite[Theorem 2.2.3]{AK}
shows that $T_p$ is isomorphic to $T_p (\Mb )$.
Thus every almost disjoint family of infinite subsets of $\N$
witnesses that $T_p$ has property $(\#)$.
Now for $1<p<\infty$, $T_p$ is complemented in $C_p$
via the projection that zeroes out the entries
that lie above the diagonal \cite{macaev}, \cite{gmf},
from which it follows easily \cite{al}
that $T_p$ is isomorphic to $C_p$.
We record these observations in Proposition \ref{Tp}.
\begin{proposition}\label{Tp}
For $1\le p \le \infty$, the space $T_p$ has property $(\#)$.
Moreover, for $1<p<\infty$, the space $C_p$ has property $(\#)$.
\end{proposition}
As we mentioned above,
it can be proved that $C_1$ and $C_{\infty}$ have the SHAI property
even though neither has an unconditional FDD.
However, the $C_p$ norms for $1\le p\le\infty$
are  what Kwapie\'n and Pe\l czy\'nski \cite{kw-p} call
unconditional matrix norms;
i.e., the norm $\left\| \sum_{i,j} a_{i,j} e_{i,j} \right\|$
of a linear combination $\sum_{i,j} a_{i,j} e_{i,j}$
of the natural basis elements
$(e_{i, j} )_{i, j = 1}^{\infty}$,
is equivalent (in our case even equal)
to the norm of $\sum_{i,j} \varepsilon_i \delta_j a_{i,j} e_{i,j}$
for all sequences of signs $(\varepsilon_i)_{i = 1}^{\infty}$
and $(\delta_j)_{j = 1}^{\infty}$.
One can define a variation of property $(\#)$
for  bases with this unconditionality property,
check that the natural bases for $C_p$, for $1\le p\le\infty$,
satisfy this property, and prove a version of Theorem \ref{SHAI}.
This shows that $C_1$ has the SHAI property
(and gives an alternative proof also for $C_p$ for $1 < p < \infty$).
This variation of Theorem \ref{SHAI} does not apply to $C_{\infty}$,
which does not have finite cotype,
and we do not know whether $C_{\infty}$ has the SHAI property.
Since our focus in this paper is on spaces
that are more closely related to $L^p$ than are the $C_p$ spaces,
we do not go into more detail.
Our main reason for bringing up $C_p$
is to point out why the definition of property $(\#)$
is made for unconditional FDDs rather than just for unconditional bases.



\noindent W.~B.~Johnson\newline
             Department Mathematics\newline
             Texas A\&M University\newline
             College Station TX 77843--3368, USA\newline
             E-mail: johnson@math.tamu.edu

  \smallskip
   \noindent N.~C.~Phillips
   \newline Department of Mathematics
   \newline University  of Oregon
   \newline Eugene OR 97403-1222, USA

\smallskip

\noindent G.~Schechtman
\newline Department of Mathematics
\newline Weizmann Institute of Science
\newline Rehovot, Israel
\newline E-mail: gideon@weizmann.ac.il


\begin{thebibliography}{33}


\bibitem{BRS}
J.~Bourgain, H.~P.~Rosenthal, and  G.~Schechtman,
\emph{An ordinal $L_p$-index for Banach spaces, with application
 to complemented subspaces of $L_p$},
Ann. of Math. (2) 114  no. 2   (1981), 193--228.

\bibitem{AK}
F.~Albiac and N.~J.~Kalton,
\emph{Topics in Banach space theory},
Graduate Texts in Mathematics, 233 (2nd edition).
Springer, New York, 2016.

\bibitem{al}
J.~Arazy and J.~Lindenstrauss,
\emph{Some linear topological properties of the spaces $C_p$
of operators on Hilbert space,}
Compositio Math. 30 (1975), 81--111.

\bibitem{barnes}
B.~A.~Barnes,
\emph{Algebraic elements of a Banach algebra modulo an ideal},
Pacific J. Math. 117, No. 2 (1985), 219--231.

\bibitem{dales}
H.~G.~Dales,
\emph{Banach Algebras and Automatic Continuity},
Oxford University Press
Inc., New York, 2000.

\bibitem{Eidh}
M.~Eidelheit,
{\emph{On isomorphisms of rings of linear operators}},
Studia Math.\  {\textbf{9}} (1940),  97--105.

\bibitem{GG}
J.~L.~B.~Gamlen and R.~J.~Gaudet,
{\emph{On subsequences of the Haar system in $L^p [-1,1]$,
$(1\le p \le \infty)$},}
Israel J. Math. 15 (1973), 404--413.

\bibitem{gmf}
I.~C.~Gohberg, A.~S.~Markus, and I.~A.~Feldman,
\emph{Normally solvable operators and
ideals associated with them,}
 American Math. Soc. Translat. 61 (1967), 63--84.


\bibitem{Horv}
B.~Horvath,
\emph{When are full representations of algebras of operators
on Banach spaces automatically faithful?},
Studia Math. 253  no. 3 (2020), 259--282.

\bibitem{HorvKan}
B.~Horvath and T.~Kania,
\emph{Surjective operators from algebras of operators
on long sequence spaces are automatically injective},
arXiv:2007.14112v1 (2020).

\bibitem{kp}
M.~I.~Kadec and   A.~Pe\l czy\'nski,
\emph{ Bases,
lacunary sequences and complemented subspaces in the spaces
 $L\sb{p}$,}
 Studia Math.  21  (1961/1962), 161--176.

\bibitem{kw-p}
S.~Kwapie\'n and   A.~Pe\l czy\'nski,
\emph{The main triangle projection in matrix spaces
and its applications,}
Studia Math. 34 (1970), 43--68.

\bibitem{LT} J.~Lindenstrauss and L.~Tzafriri,
\emph{Classical Banach Spaces, Vol I }, Ergebn. Math.
Grenzgeb. {92}, Springer, (1977).

\bibitem{macaev}
V.~I.~Macaev,
\emph{Volterra operators obtained from self-adjoint operators
by perturbation,}
(Russian) Dokl. Akad. Nauk SSSR 139 (1961), 810--813.

\bibitem{ro}
H.~P.~Rosenthal,
\emph{On the subspaces of $L^p$ $(p>2)$ spanned by sequences
of independent random variables,}
Israel J. Math.  8  (1970), 273--303.

\bibitem{sc}
G.~Schechtman,
\emph{ Examples of $\mathcal{L}_p$ spaces $(1<p\neq 2<\infty)$,}
Israel J. Math.  22  (1975),  no. 2, 138--147.

\end{thebibliography}
\end{document}